\documentclass{amsart}

\usepackage{amssymb}
\usepackage[cp1251]{inputenc}
\usepackage[english]{babel}
\usepackage{amsmath}
\usepackage{color}
\usepackage{graphics}
\usepackage{epsfig}
\usepackage[all,cmtip]{xy}
\usepackage{pb-diagram,pb-xy}
\usepackage{bm}
\usepackage{multicol}

\newtheorem{thm}{Theorem}[section]

 \newtheorem{proposition}[thm]{Proposition}
 \theoremstyle{definition}
 \newtheorem{definition}[thm]{Definition}
 \theoremstyle{remark}

 \numberwithin{equation}{section}
 
\begin{document}

%-------------------------------------------------------------------------
% editorial commands: to be inserted by the editorial office
%
%\firstpage{1} \volume{228} \Copyrightyear{2004} \DOI{003-0001}
%\seriesextraline{This is the Concrete Title of this Book\br H.E. R and S.T.C. W, Eds.}
%\firstpage{1}

\title[Orbit groups]
 {Orbit groups}

%----------Author 1
\author[Quitzeh Morales Mel\'endez]{Quitzeh Morales Mel\'endez}

\address{%
Camino a la Zanjita s/n\\
Noche Buena\\
Sta. Cruz Xoxocotl\'an\\
71230 Oaxaca\\
M\'exico}

\email{qmoralesme@conacyt.mx, qmorales@upn.mx}

\thanks{This work was completed with the support of C\'atedras CONACYT project 1522.}
%----------Author 2
%\author{A Second Author}
%\address{The address of\br
%the second author\br
%sitting somewhere\br
%in the world}
%\email{dont@know.who.knows}


\date{May 13, 2020}
%\dedicatory{To Nayar}
%%% ----------------------------------------------------------------------

\begin{abstract}
In the paper \cite{AJMM} are introduced two groups generated by the orbits of an action of a group on another group by automorphisms. One is of group-theoretic nature and the other comes from homology of invariant group chains. In this note are  given some properties of the first groups and is studied a natural homomorphism between these groups. More precisely, it is shown that this homomorphism is not injective nor surjective. A description of the kernel is given.
\end{abstract}

%%% ----------------------------------------------------------------------
\maketitle
%%% ----------------------------------------------------------------------
%\tableofcontents
\section*{Introduction}
The study of orbits of group actions on groups by group automorphism has been related to the study of finite solvable groups, in particular that of Frobenius groups 
(See \cite{Thompson59}, \cite{Feit}). More recent work of Deaconescu and Walls is dedicated to the study of orbits and fixed points of such actions \cite{DeaconescuWalls}.

The present work has another origin coming from the study of homology theories for groups and their automorphism groups. Namely, it was found in \cite[Theorem 2]{JLM} a group generated by the orbits of the action of a finite group on another group. This group arrises as the first homology group of invariant group chains defined by Knudson \cite{Knudson} and, as expected, it is related to the generalization in this context of the abelianization of a group \cite[Sec. 2]{AJMM}.

The discovery of an abelian group generated by orbits of these actions rises the question on the existence of such an object in group theory. This is defined in \cite[Sec. 2.2]{AJMM} together with a homomorphism to the first homology group of invariant group chains.

The first part of this work is dedicated to discuss some properties of the so-called orbit groups and the second part to study relations with the first homology group of invariant group chains of Knudson. It is shown that the homomorphism defined in \cite[Sec. 2.2]{AJMM} is not in general surjective nor injective. Condition for elements in the kernel of this homomorphism are given.

\section{Orbit groups and main properties}
Consider the category having as objets groups with a group action by group automorphisms and as morphisms pairs of homomorphisms commuting with the actions. More precisely, the objects of this category are triples $(Q, G, \varphi)$ where $Q$ and $G$ are groups and 
$\varphi: Q\times G \to G$ is a group action by group homomorphisms, i.e.
$\varphi(q_1q_2,g)=\varphi(q_1,\varphi(q_2,g)), \varphi(1,g)=g$ and 
$\varphi(q,g_1g_2)=\varphi(q,g_1)\varphi(q,g_2), \varphi(q,1)=1$, or equivalently,
$$
(q_1q_2)(g)=q_1(q_2(g)),  \qquad 1(g)=g,
$$and 
$$
 q(g_1g_2)=q(g_1)q(g_2) \qquad q(1)= 1
$$ for $q, q_1, q_2\in Q$, $g, g_1, g_2\in G$, if the action is avoided from the notation.
The morphisms of this category are pairs $(\alpha, \beta )$ of homomorphisms 
$\alpha: Q\to Q', \,\beta: G\to G'$ commuting with the actions, i.e. such that the diagram
$$
\xymatrix{Q\times G   \ar[r]^-{\varphi} \ar[d]_-{\alpha \times \beta} & G\ar[d]_-{\beta}\\ 
                Q'\times G' \ar[r]^-{\varphi'} &  G'} 
$$is commutative.

In this case, one has the well known semidirect product functor $Q\ltimes_\varphi G$ (or $Q\ltimes G$)  to the category of groups. In order to keep track of the group $G$, one can consider the corresponding functor to the category of pairs of groups. This functor, say $L$, takes the triple $(Q, G, \varphi)$ to the pair $(G,Q\ltimes G)$ and the pair of homomorphisms $(\alpha, \beta )$ to the pair 
$(\beta, \alpha \ltimes \beta)$.

\begin{definition}[Conf. \cite{AJMM}]
The orbit group of the action of $Q$ on $G$, denoted by $G//Q$, is a group equipped with an surjective homomorphism $p: G \to G//Q$ and having the property that every group homomorphism $\phi: G \to H$ identifying orbits, i.e. such that $\phi(q(g))=\phi(g), g\in G$, factors trough $G//Q$, i.e. there exists an unique homomorphism $\psi:G//Q \to H$ such that the diagram 
$$
\xymatrix{G  \ar[dr]^-{\varphi} \ar[d]_-{p} & \\ 
                G//Q \ar[r]^-{\psi} &    H,} 
$$is commutative.
\end{definition}

In order to show the existence of such an object, construct first a quotient of the semidirect product $Q\ltimes G$ by a normal subgroup, say $N$, such that the relation $\bar g=\overline{q(g)}, g\in G, q\in Q$ is true in $(Q\ltimes G)/N$. Such subgroup $N$ must contain the elements $gq(g^{-1})$ for every $ g\in G, q\in Q$. In the semidirect product, the action of $Q$ on $G$ can be written as conjugation by the corresponding element. In this way, one has
$gq(g^{-1})=[g,q]$ and, therefore, the commutator subgroup $[G,Q]$ must be contained in $N$. However, the subgroup $[G,Q]$ is not necessarily normal, so, one shall take the minimal normal subgroup containing $[G,Q]$, this is $N=[G,Q]^{Q \ltimes G}$.  This gives the quotient 
$Q \ltimes G/[G,Q]^{Q \ltimes G}$ together with the quotient map  
$ Q \ltimes G \longrightarrow (Q \ltimes G)/[G,Q]^{Q \ltimes G}$. The map $p: G \to G//Q$ is then obtained as the restriction of the composition 
$$
G\hookrightarrow Q \ltimes G \longrightarrow (Q \ltimes G)/[G,Q]^{Q \ltimes G}
$$
to its image. By the second isomorphism theorem, one obtains the formula 
$G//Q=G/(G\cap [G,Q]^{Q \ltimes G})$. Using the exterior definition of the semidirect product, one can show that $[G,Q]^{Q \ltimes G}$ is contained in $G$ and coincides with $[G,Q]^{G}$. Therefore, one obtains a simpler expression $G//Q=G/[G,Q]^{G}$ which is defined solely in terms of $G$ and the action of $Q$ on $G$, because the group  
$[G,Q]^{G}$ is generated by the elements of the form $g_1gq(g^{-1})g_1^{-1}$ for $ g, g_1\in G, q\in Q$.

The universal property of $G//Q$ means that a morphism 
$$(\alpha, \beta): (Q, G, \varphi)\longrightarrow (Q', G', \varphi')$$ induces a homomorphism 
$\beta//\alpha: G//Q \rightarrow G'//Q'  $ such that the square 
$$
\xymatrix{G  \ar[r]^-{\beta} \ar[d]_-{p} & G'\ar[d]_-{p'}\\ 
                G//Q \ar[r]^-{\beta//\alpha} &  G'//Q',} 
$$is commutative. In other words, one has a functor from the category of groups acting on groups to the category of groups.

The orbit group functor just constructed has the following properties:
\begin{itemize}
\item if $i:K\to G$ is a monomorphism and $K$ is a $Q$-invariant subgroup, then $i//id:K//Q\to G//Q$ is a monomorphism;
\item if $j:G\to H$ is an equivariant epimorphism, then $j//id:H//Q\to G//Q$ is an epimorphism;
\item if $N \lhd G $ is a $Q$-invariant normal subgroup, then $N//Q \lhd G//Q$ is a normal subgroup and one has
$
(G/N)//Q\cong (G//Q)/(N//Q);
$
\item there is a correspondence between subgroups of $H'< G//Q$ and $Q$-invariant subgroups 
 $H< G$ containing $[G,Q]^{G}$;
 \item there is a correspondence between quotients of $(G//Q)/H'$ and quotients $G/H $ by $Q$-invariant subgroups such the induced action $Q\times G/H\to G/H$ is trivial.
\end{itemize}

There is no natural transformation between the functors $Q \ltimes G$ and $G//Q$ in the very manner they are defined. Instead, one can consider the corresponding functors $(G, Q \ltimes G) $ and  $(G//Q,(Q \ltimes G)/[G,Q]^{G})$ to the category of pairs of groups. Then, the pair of morphisms $(p,r) : (G, Q \ltimes G) \rightarrow (G//Q,(Q \ltimes G)/[G,Q]^{G})$, where 
$r :Q \ltimes G \rightarrow (Q \ltimes G)/[G,Q]^{G}$ is the quotient map, is natural.

The functor $G//Q$ is a (non abelian) generalization of the abelianization of a group: take $Q=G$ acting on itself by inner automorphisms, then one has $G//G=G/[G,G]=G_{\text{ab}}$.

\section{Relations between the orbit group and the first homology group of invariant group chains}
Let $Q$ be a finite group, $G$ be a group, and $Q\times G \to G $ be an action of $Q$ on $G$ by 
group automorphisms. Denote by $H_{1}^{Q}(G,\mathbb{Z})$ the first homology group of $Q$-invariant chains on the group $G$ (see \cite{Knudson}). It was shown in \cite[Theorem 2]{JLM} the following.
\begin{eqnarray*}
&&H^{Q}_{1}(G,\mathbb{Z})=\\
&&\frac{\mathbb{Z}\left\lbrace
\sum_{\bar{q}\in Q/Q_{[g]}}q[g]\,\vline\, g\in G  \right\rbrace}{\mathbb{Z}\left\lbrace
a\sum_{\bar{q}\in Q/Q_{[g_{2}]}}q[g_{2}]-
b\sum_{\bar{q}\in Q/Q_{[g_{1}g_{2}]}}q[g_{1}g_{2}]
+c\sum_{\bar{q}\in Q/Q_{[g_{1}]}}q[g_{1}]  \right\rbrace},
\end{eqnarray*}
where $a=\frac{\mid Q_{[g_{2}]}\mid}{\mid Q_{[g_{1}]}\cap Q_{[g_{2}]}\mid}$, $b=\frac{\mid Q_{[g_{1}g_{2}]}\mid}{\mid Q_{[g_{1}]}\cap Q_{[g_{2}]}\mid}$,$c=\frac{\mid Q_{[g_{1}]}\mid}{\mid Q_{[g_{1}]}\cap Q_{[g_{2}]}\mid}$, $g_1, g_2\in G$.

This result means that the abelian group $H^{Q}_{1}(G,\mathbb{Z})$ is generated by 
the orbits of the action of $Q$ on $G$.

Let $C_1(G)$ be the free abelian group generated by $G$ as a set. This group has a natural action of $Q$ and the norm map $N: G\to C_1(G)^Q$ sending each element $g\in G$ to the sum of the elements on its orbit  $$\sum_{[q]\in Q}q[g]=\mid Q_{g}\mid\sum_{[q]\in Q/Q_{g}}q[g]$$ identifies $Q$-orbits and, by the universal property of the orbit group and the abelianization of a group, this induces a homomorphism (conf.\cite{AJMM})
\begin{equation}\label{homomorphism}
(G//Q)_{\text{ab}}\longrightarrow H^{Q}_{1}(G,\mathbb{Z})
\end{equation} such that the diagram
$$\xymatrix{G  \ar[ddr]^-{\bar N} \ar[d]_-{p} & \\ 
                G//Q \ar[d] &   \\
                 (G//Q)_{\text{ab}} \ar[r] & H^{Q}_{1}(G,\mathbb{Z} ),}$$
is commutative.

It is intuitively clear that homomorphism \eqref{homomorphism} may not in general be surjective, as it is induced by the non surjective norm map. 

In order to have concrete examples of this, consider the case $G=\mathbb{Z}/4$, 
$Q=\mathbb{Z}/2=\langle t \rangle$ acting by group inversion, i.e. $n \mapsto - n$. In this case, 
$(G//Q)_{\text{ab}}=G//Q=G/[G,Q]^G=G/[G,Q]$ because $G$ is abelian. 
The subgroup $[G,Q]$ is generated by elements of the form $n + t(-n)=2n$, i.e.  
$[G,Q]=\{0,2\}$ and $(G//Q)_{\text{ab}}\cong \mathbb{Z}/2$. On the other hand, according to \cite[Theorem 3]{JLM}, $H^{Q}_{1}(G,\mathbb{Z} )\cong \mathbb{Z}/2\oplus \mathbb{Z}/2$ and, therefore, homomorphism \eqref{homomorphism} can not be surjective.

%It might be not injective also: using the relation (\ref{minimalisotropyrelation}), it is easy to show that the order $n$ of an element $g\in G$ annihilates the corresponding element $[g]^Q\in  H^{Q}_{1}(G,\mathbb{Z} )$. Therefore, if $n$ divides  $\mid Q_g\mid$, then the image of such element would be zero. 

\section{Conditions for (non) inyectivity of homomorphism \eqref{homomorphism} }

The aim of this section is to give conditions for the existence of elements in the kernel of homomorphism \eqref{homomorphism}
%$(G//Q)_{\text{ab}}\longrightarrow H^{Q}_{1}(G,\mathbb{Z})$ 
and to give more detailed conditions for these elements when the order of the group $Q$ is bounded. 
%In particular, one can show that this homomorphism is inyective for Frobenius groups.

Assume that the element $g\in G$ is in the kernel of the map $\tilde N : G \to H^{Q}_{1}(G,\mathbb{Z})$. This means that, in the free group $C_1(G)^Q$, one has an equation of the form
\begin{equation*}
\sum_{q\in Q}q(g)=
\end{equation*}
\begin{equation}\label{positivecoefficients}
\sum_{i=1}^{N}\left[ a_i\sum_{q\in Q/Q_{g_{1i}}}q(g_{1i})+
b_i\sum_{q\in Q/Q_{g_{2i}}}q(g_{2i})-c_i\sum_{q\in Q/Q_{g_{1i}g_{2i}}}q(g_{1i}g_{2i}) \right] 
\end{equation}where $a_i=\pm\frac{|Q_{g_{1i}} |}{|Q_{g_{1i}}\cap Q_{g_{2i}} |}$, 
$b_i=\pm\frac{|Q_{g_{2i}} |}{|Q_{g_{1i}}\cap Q_{g_{2i}} |}$ and
$c_i=\mp\frac{|Q_{g_{1i}g_{2i}} |}{|Q_{g_{1i}}\cap Q_{g_{2i}} |}$.
Equivalently, 
\begin{multline}\label{Inkernel}
\sum_{i=1}^{N}\left[ a_i\sum_{q\in Q/Q_{g_{1i}}}q(g_{1i})+
b_i\sum_{q\in Q/Q_{g_{2i}}}q(g_{2i})-c_i\sum_{q\in Q/Q_{g_{1i}g_{2i}}}q(g_{1i}g_{2i}) \right]-\\ - \sum_{q\in Q}q(g)=0.  
\end{multline}
This equation means that all terms with negative sign must coincide with some elements with positive sign. 

In order to find elements in the kernel of homomorphism \eqref{homomorphism}, one may consider special cases of equation \eqref{Inkernel}. In the following, it will be considered that coefficients $a_i$ and $b_i$ are all positive and $c_i$ are all negative.

A necessary condition for equation  \eqref{Inkernel} to hold is that the number of positive terms must coincide with the number of negative terms. In this case,
\begin{equation}
\sum_{i=1}^{N}( a_i [Q:Q_{g_{1i}}] +b_i [Q:Q_{g_{2i}}])=  |Q| +
\sum_{i=1}^{N} c_i[Q:Q_{g_{1i}g_{2i}}] .  
\end{equation}Equivalently, one has
\begin{multline}
\sum_{i=1}^{N}\left( \frac{|Q_{g_{1i}} |}{|Q_{g_{1i}}\cap Q_{g_{2i}} |}\frac{|Q|}{|Q_{g_{1i}} |} +\frac{|Q_{g_{2i}} |}{|Q_{g_{1i}}\cap Q_{g_{2i}} |}\frac{|Q|}{|Q_{g_{2i}} |}\right)=\\
=\sum_{i=1}^{N} \frac{|Q_{g_{1i}g_{2i}} |}{|Q_{g_{1i}}\cap Q_{g_{2i}} |}\frac{|Q|}{|Q_{g_{1i}g_{2i}} |} + |Q|,  
\end{multline}i.e.
\begin{equation}
\sum_{i=1}^{N}\left( \frac{|Q |}{|Q_{g_{1i}}\cap Q_{g_{2i}} |} +\frac{|Q|}{|Q_{g_{1i}}\cap Q_{g_{2i}} |}\right)=
\sum_{i=1}^{N} \frac{|Q|}{|Q_{g_{1i}}\cap Q_{g_{2i}} |} + |Q|,  
\end{equation}that is 
\begin{equation}
\sum_{i=1}^{N}\left( 2\frac{|Q |}{|Q_{g_{1i}}\cap Q_{g_{2i}} |} \right)=
\sum_{i=1}^{N} \frac{|Q|}{|Q_{g_{1i}}\cap Q_{g_{2i}} |} + |Q|,  
\end{equation}
\begin{equation}
2\sum_{i=1}^{N}\left( \frac{|Q |}{|Q_{g_{1i}}\cap Q_{g_{2i}} |} \right)=
\sum_{i=1}^{N} \frac{|Q|}{|Q_{g_{1i}}\cap Q_{g_{2i}} |} + |Q|,  
\end{equation}
\begin{equation}
\sum_{i=1}^{N}\left( \frac{|Q |}{|Q_{g_{1i}}\cap Q_{g_{2i}} |} \right)=|Q|.
\end{equation}Therefore, one has the equation

\begin{equation}\label{partitionof1}
\sum_{i=1}^{N}\left( \frac{1}{|Q_{g_{1i}}\cap Q_{g_{2i}} |} \right)=1.
\end{equation}That is, intersection of isotropy groups of elements in right side of equation \eqref{Inkernel} must form a rational partition of 1. 

Equation \eqref{partitionof1} can be applied to bound the length $N$ of the sum in equation 
\eqref{Inkernel}. For example, if $Q$ does not fix elements different from the identity in $G$, then one has $Q_{g_{1}}\cap Q_{g_{2}}=1$ for every pair of elements $g_1, g_2\in G$ that are not the neutral element. In this case one must have $N=1$, i.e. equation \eqref{Inkernel} reduces to 
\begin{equation}\label{length1}
0=  a\sum_{q\in Q/Q_{g_{1}}}q(g_{1})+
b\sum_{q\in Q/Q_{g_{2}}}q(g_{2})-c\sum_{q\in Q/Q_{g_{1}g_{2}}}q(g_{1}g_{2}) -\sum_{q\in Q}q(g).  
\end{equation}

For the equation to hold, $g_1$ or $g_2$ must be in the same orbit with their product $g_1g_2$,  for example, one has $q_1(g_2)=g_1g_2$, which, on the one hand, implies that  $|Q_{g_{2}}|=|Q_{g_1g_{2}}|$, i.e.
$b=c$, and on the other hand, that $\sum_{q\in Q/Q_{g_1}}q(g_1)=\sum_{q\in Q/Q_{g_1g_2}}q(g_1g_2)$. Therefore, equation \eqref{length1} becomes 
\begin{equation}
0=  a\sum_{q\in Q/Q_{g_{1}}}q(g_{1})-\sum_{q\in Q}q(g).
\end{equation}This means that $g$ and $g_1$ are in the same orbit, but one has that
$g_1=g_2^{-1}q_1(g_2)=[g_2,q_1]^{-1}$, i.e. $\bar g=\bar g_1\in G//Q$. And this does not give non-trivial elements in the kernel.

%Therefore homomorphism (\ref{homomorphism}) is injective for free actions, as in the case of Frobenius groups (see, for example, \cite{Thompson59}).

The length $N$ may also be bounded in terms of the order of the group, for example, if $|Q|=2$, then the only possible non trivial length (in terms of the previous discussion) is $N=2$. The only corresponding partition in this case is $\frac{1}{2}+\frac{1}{2}$.

\begin{proposition}
Let $Q$ be a group of order two. If a non trivial element $\bar g \in (G//Q)_{\text{ab}}$ satisfies equation \eqref{positivecoefficients} then it must have order two.
\end{proposition}
\begin{proof}
As it was shown, for length $N=1$ one has  $\bar g=1 \in G//Q$ which implies $\bar g=1 \in (G//Q)_{\text{ab}}$. So, the only remaining possibility is $N=2$.

Thus, one has the equation
\begin{equation}
\begin{array}{c}\label{length2}
a_1\sum_{q\in Q/Q_{g_{11}}}q(g_{11})+
b_1\sum_{q\in Q/Q_{g_{21}}}q(g_{21})-c_1\sum_{q\in Q/Q_{g_{11}g_{21}}}q(g_{11}g_{21})+ \\
a_2\sum_{q\in Q/Q_{g_{12}}}q(g_{12})+
b_2\sum_{q\in Q/Q_{g_{22}}}q(g_{22})-c_2\sum_{q\in Q/Q_{g_{12}g_{22}}}q(g_{12}g_{22}) \\
-\sum_{q\in Q}q(g)=0. 
\end{array}
\end{equation}

\textbf{Case 1.}  \textit{A cancellation in the upper sum in equation \eqref{length2}.} As in the previous 
discussion, in this case one has $b_1=c_1$, the corresponding terms in the upper sum cancel, $\bar g_{11}=1 \in G//Q$ and one obtains the equation
\begin{equation}
\begin{array}{c}
0= a_1\sum_{q\in Q/Q_{g_{11}}}q(g_{11})+ \\
a_2\sum_{q\in Q/Q_{g_{12}}}q(g_{12})+
b_2\sum_{q\in Q/Q_{g_{22}}}q(g_{22})-c_2\sum_{q\in Q/Q_{g_{12}g_{22}}}q(g_{12}g_{22}) \\
-\sum_{q\in Q}q(g). 
\end{array}
\end{equation} If the lower therm cancels with (part of) the upper remaining term, then $g$ and 
$g_{11}$ are and the same orbit and $\bar g=1 \in (G//Q)_{\text{ab}}$. Otherwise, the upper term cancels with (part of) the negative middle term which implies $|Q_{g_{11}}|=|Q_{g_{12}g_{22}}|$ and, as $|Q_{g_{11}}\cap Q_{g_{21}} |=2=|Q_{g_{12}}\cap Q_{g_{22}} |$ one has $a_1=c_2$ and 
both remaining middle (positive) terms must cancel with the lower therm. This means that both $g_{12}$ and $g_{22}$ are in the same orbit of $g$ and, therefore
$$ (\bar g)^2=\bar g\bar g= \bar g_{12} \bar g_{22}=\overline{g_{12}  g_{22}}= \bar g_{11} =1\in  G//Q. $$ 
\textbf{Case 2.}  \textit{No cancellations in the same row in equation \eqref{length2}.} One may assume that there exist elements $q_1, q_2 \in Q$ such that $q_1(g_{11})=g_{12}g_{22} $ and 
$q_2(g_{12})=g_{11}g_{21} $. Then 
$$
q_2q_1(g_{11})= q_2(g_{12}g_{22})=q_2(g_{12})q_2(g_{22})=g_{11}g_{21}q_2(g_{22}),
$$which implies $[g_{11}^{-1}, q_2q_1]=g_{11}^{-1}q_2q_1(g_{11})=g_{21}q_2(g_{22})$, i.e.
$\bar g_{21}\bar g_{22}=1\in G//Q$.

As in case 1, one has $a_1=c_2$ and analogously $a_2=c_1$. The remaining positive terms must cancel with $-\sum_{q\in Q}q(g)$ which implies that $g_{21}$, $ g_{22}$ and $g$ are in the same orbit. So, $\bar g \in G//Q$ is trivial or it has order two.
\end{proof}

Note that, if $|Q|=2$, then $|Q_{g_{11}}\cap Q_{g_{21}} |=2=|Q_{g_{12}}\cap Q_{g_{22}} |$ implies 
$Q_{g_{12}}\cap Q_{g_{22}}=Q$ and, so, $Q_{g_{12}}= Q_{g_{22}}=Q$. Similarly, all the other isotropy groups involved in equation \eqref{length2} are equal to $Q$. This means that elements in the kernel of \eqref{homomorphism} from equation \eqref{positivecoefficients} must come from order two elements in the fixed point subgroup $G^Q$ ($=C_Q(G)$ in the notation used in \cite{DeaconescuWalls}) of $G$.

As an example of this, consider $G=\mathbb{Z}/4$ with $Q$ acting as inversion 
$n\mapsto -n$. Then the equation,
$$
0=([2]+[2] - [0]) + ([0]+[2] - [2]) -2 [2]
$$
shows that the class of the (non trivial) element $g=2$ in $(G//Q)_{\text{ab}}$ is in the kernel of \eqref{homomorphism}. This example can be generalized as follows.
\begin{proposition}\label{fixedordern}
If $|Q|=n$ and there is a fixed element $g\in G$ of order $n$, then the class of this element in $(G//Q)_{\text{ab}}$  is in the kernel of \eqref{homomorphism}.
\end{proposition}
\begin{proof}
As before, one has the length $n$ equation
$$
0=([g]+    [g] - [g^2]) + 
    ([g^2]+[g] - [g^3]) 
    + \cdots $$ 
    $$\cdots+
    ([g^{n-1}]+[g] - [1])
    + ([1]+[g] - [g]) - n [g].
$$\end{proof}

Before returning to the general case of equation \eqref{positivecoefficients}, note that 
the class $[1]$ of the neutral element in $H^{Q}_{1}(G,\mathbb{Z})$ must be zero. This means that one may consider equation \eqref{positivecoefficients} in the quotient $C^{Q}_{1}(G)/\mathbb{Z}[1]$ instead. 

Considering this, if some $g_{1i}=1$, then  $Q_{g_{1i}}=Q$ and $Q_{g_{1i}}\cap Q_{g_{2i}}=Q_{g_{2i}}$ giving $b_i=c_i=1$. So, the corresponding term of the sum is zero:
$$ a_i\sum_{q\in Q/Q_{g_{1i}}}q(g_{1i})+
b_i\sum_{q\in Q/Q_{g_{2i}}}q(g_{2i})-c_i\sum_{q\in Q/Q_{g_{1i}g_{2i}}}q(g_{1i}g_{2i})
=$$
$$ 0 + \sum_{q\in Q/Q_{g_{2i}}}q(g_{2i})-\sum_{q\in Q/Q_{1\cdot g_{2i}}}q(1\cdot g_{2i})=0,$$
and one may assume $g_{1i}\not = 1 \not = g_{2i}$. 

Also, if $g_{2i}=g_{1i}^{-1}$ then these elements have the same isotropy group, so $a_i=b_i=1$ and the corresponding term of the sum is just 
$ [g_{1i}] + [g_{1i}^{-1}]$. 

Now, in equation \eqref{positivecoefficients} one can rearrange the sum in the right side to have
\begin{multline}\label{allcoefficients}
\sum_{q\in Q}q(g) = \\
=\sum_{i=1}^{N}\left[ a_i\sum_{q\in Q/Q_{g_{1i}}}q(g_{1i})+b_i\sum_{q\in Q/Q_{g_{2i}}}q(g_{2i})-c_i\sum_{q\in Q/Q_{g_{1i}g_{2i}}}q(g_{1i}g_{2i}) \right] \\ - \sum_{j=1}^{M}\left[ d_i\sum_{q\in Q/Q_{h_{1i}}}q(h_{1i})+
e_i\sum_{q\in Q/Q_{h_{2i}}}q(h_{2i})-f_i\sum_{q\in Q/Q_{h_{1i}h_{2i}}}q(h_{1i}h_{2i}) \right] 
\end{multline}
Counting the number of elements in the sum, one obtains the equation
\begin{equation}\label{difference_partitionof1}
1= \sum_{i=1}^{N}\left( \frac{1}{|Q_{g_{1i}}\cap Q_{g_{2i}} |} \right)
-\sum_{j=1}^{M}\left( \frac{1}{|Q_{h_{1i}}\cap Q_{h_{2i}} |} \right).
\end{equation}where all elements $g_{1i}, g_{2i}, h_{1i}, h_{2i}\in G$ are different
from the neutral element. For free actions, for example, one obtains the equation $1=N-M$.

Consider the set-theoretic quotient map $p:G\to G/Q$, where the elements of this set are represented by sums of the form $x=\sum_{[q]\in Q/Q_{g'}}q(g)$, then all elements in the same inverse image of $p$ can be summed up. This gives an equation of the form
$$
 [\sum_{i, p(g_{1i})=p(g)} a_i+\sum_{i, p(g_{2i})=p(g)} b_i-
\sum_{i, p(g_{1i}g_{2i})=p(g)} c_i -  
$$
$$
\sum_{i, p(h_{1i})=p(g)}d_i -\sum_{i, p(h_{2i})=p(g)}e_i+\sum_{i, p(h_{1i}h_{2i})=p(g)}f_i -\mid Q_g\mid] [g]+
$$
$$
+\sum_{x\in G/Q} [\sum_{i, p(g_{1i})=x} a_i+\sum_{i, p(g_{2i})=x} b_i-
\sum_{i, p(g_{1i}g_{2i})=x} c_i -  
$$
$$
\sum_{i, p(h_{1i})=x}d_i -\sum_{i, p(h_{2i})=x}e_i+\sum_{i, p(h_{1i}h_{2i})=x}f_i ]x=0.
$$As the generators of the free group $C^Q_1(G)$ are independent, this equation can only be true if the coefficients are all zero, giving the equations:
\begin{multline}
 \sum_{i, p(g_{1i})=p(g)} a_i+\sum_{i, p(g_{2i})=p(g)} b_i-
\sum_{i, p(g_{1i}g_{2i})=p(g)} c_i -  
\\
\sum_{i, p(h_{1i})=p(g)}d_i -\sum_{i, p(h_{2i})=p(g)}e_i+\sum_{i, p(h_{1i}h_{2i})=p(g)}f_i -\mid Q_g\mid = 0;
\\
\sum_{i, p(g_{1i})=x} a_i+\sum_{i, p(g_{2i})=x} b_i-
\sum_{i, p(g_{1i}g_{2i})=x} c_i -  
\\
\sum_{i, p(h_{1i})=x}d_i -\sum_{i, p(h_{2i})=x}e_i+\sum_{i, p(h_{1i}h_{2i})=x}f_i =0.
\end{multline} where all elements $g_{1i}, g_{2i}, h_{1i}, h_{2i}\in G$ are different
from the neutral element. These equations, together with condition 
\eqref{allcoefficients} give a description of the kernel of \eqref{homomorphism}.

%\begin{rem}
%Additional comments are being typeset without boldfaced entrance
%word as they may be minor important.
%\end{rem}

%\begin{ex}
%For some constructs, even no number is required.
%\end{ex}

%Displayed equations may be numbered like the following one:
%\begin{equation}
%\sqrt{1-\sin^2(x)}=|\cos(x)|.
%\end{equation}

%\paragraph{Here is a Sample for a Paragraph}

%As you observe, paragraphs do not have numbers and start new lines after the heading, by default.
% ------------------------------------------------------------------------

%\subsection*{Acknowledgment}
%Many thanks to our \TeX-pert for developing this class file.

% ------------------------------------------------------------------------
\end{document}